%% file: q4.tex
\documentclass[a4paper,twoside,reqno]{amsart}
\calclayout
\usepackage{hyperref}
\usepackage{url}
\usepackage{amssymb,amscd}
\usepackage[mathscr]{eucal}
\usepackage{thmtools}
\usepackage{cleveref}
\usepackage[section]{placeins}
\usepackage{amsaddr}
\input lst-Macaulay2.tex
\date{\today}
\DeclareMathOperator{\Tr}{Tr}
\DeclareMathOperator{\diag}{diag}

\newcommand{\RR}{\mathbb{R}}

\newcommand{\HF}[3]{H_{{#2},{#1}}({#3})}
\newcommand{\PHF}[2]{P_{{#2},{#1}}}
\newcommand{\SosHF}[2]{\Sigma_{{#2},{#1}}}

\newcommand{\PH}[2]{P_{{#1},{#2}}}
\newcommand{\SosH}[2]{\Sigma_{{#1},{#2}}}

\theoremstyle{plain}
\newtheorem{prop}{Proposition}
\newtheorem{thm}[prop]{Theorem}

\newtheorem{lem}[prop]{Lemma}

\theoremstyle{definition}

\theoremstyle{remark}

\newtheorem{Examples}[prop]{Examples}
\newtheorem{rem}[prop]{Remark}
\newtheorem{Rems}[prop]{Remarks}

\newcounter{reml}

\newenvironment{remlist}{\begingroup\setcounter{reml}{0}}
  {\endgroup}

\title[Efficient sum of squares certificate for 4-ary 4-ic]
{An efficient sum of squares nonnegativity certificate for quaternary quartic}
\author{Dmitrii Pasechnik}
\address{Department of Computer Science, Northwestern University, Evanston, Illinois, USA}

\begin{document}
\begin{abstract}
For any 4-variate quartic form $f\geq 0$ (i.e. $f$ nonnegative, homogeneous polynomial of degree $4$
with real coefficients) there exist 
quadratic forms $q$ and $q'$ so that $qq'f$ is a sum of squares (s.o.s.) of quartics,
by reducing to the case of $f=au^2+2bu+c$  with $a$, $b$, $c$ $3$-variate forms of degrees 2, 3, 4, respectively,
and invoking  on its discriminant $\Delta=ac-b^2$ a theorem by Hilbert (1893) asserting that for any ternary sextic
$h\geq 0$ there exists a quadric $q''$ so that $q''h$ is s.o.s. of quartics.

Towards deciding whether
just one $q$ always suffices to make $qf$ a s.o.s, we give explicit examples of non-s.o.s.
$f=au^2+2bu+c\geq 0$ with non-s.o.s. $\Delta$. However, in all these examples $af$ are s.o.s. That is, the straightforward s.o.s. decomposition via Hilbert (1893) need not be the best possible.
While it remains open whether one $q$ always suffices (and we conjecture that $q=a$ suffices), we describe how
the existence of such $q$ is related to particular
types of s.o.s. decompositions for $\Delta$. 
\end{abstract}

\maketitle

\section{Introduction}
Certifying global nonnegativity of an $n$-variate polynomial $f$ of degree $2d$ 
($\deg{f}=2d$, $f\geq 0$, for short) by
decomposing it as a sum of squares of ``nice'' functions is a classical topic.
In 1888 Hilbert \cite{Hilb-deg4} has shown that these ``nice'' functions cannot always be
just polynomials, unless $n=1$, or $2d=2$, or $2d=4$, $n=2$ (see e.g. \cite{MR2500469} for a modern
exposition).
A less known Hilbert's paper \cite{Hilb-ternsos} from 1893 shows that for $n=2$ one can decompose $f\geq 0$ as
a sum of squares of rational functions, with degrees of denominators (and thus numerators)
bounded from above by $O(d^2)$.
In particular, \cite{Hilb-ternsos} established that for any $h\geq 0$, $\deg{h}=2d$, $n=2$ there exists $q\geq 0$, $\deg{q}=2d-4$, s.t. $qh$ is an s.o.s. of degree $2d-2$ polynomials. In view of \cite{Hilb-deg4}, 
the bound on $\deg{q}$ is sharp for $2d=6$.

\begin{rem}
Blekherman in \cite{MR3272733} provided an alternative proof of the latter sharp bound (Theorem~2.6, with a proof on p.81).
Further improvements on sharpness for $d>3$ can be found in Blekherman, Smith, Velasco \cite{MR3998790} and in
Blekherman, Sinn, Smith, Velasco \cite{blekherman2024}.

In \cite{MR3272733} (cf. Question after Corollary~2.8 there) it was asked whether for any $h\geq 0$ with $n=3$, $2d=4$ there exists a quadratic multiplier $q$ so that $qh$ is an s.o.s.
Here we show that two quadratic multipliers suffice.
\end{rem}

Tools at Hilbert's disposal were insufficient to
extend this to $n>2$; this led to the 17-th Problem from Hilbert's celebrated
list \cite{Hilb-probs}, 
settled in the affirmative by Artin and Schrieier
\cite{MR3069468}. Namely, $f\geq 0$ is decomposable
as a sum of squares of rational functions. However, first \emph{bounds} on degrees of the corresponding
denominators for $n>2$ only appeared much later; the best presently known bound is a height $5$ tower of exponents,
involving $n$ and $d$, see \cite{MR4071235}.
Here we take on the first case without a nice bound, namely, $n=3$, $d=4$, and prove
an almost sharp bound, by surprisingly elementary means. 

Technically, it is more natural to work with \emph{forms}, i.e. homogeneous polynomials,
as homogenisation does not change the s.o.s. decomposability.
Consider $\RR[x_1,\dots,x_n]:=\RR[x]$ and its homogeneous subspaces $\HF{d}{n}{\RR}$
of degree $d$. Then the nonnegative on $\RR^{n}$ polynomials $f\in \HF{2d}{n}{\RR}$
form a convex cone, denoted by $\PHF{2d}{n}$.
A $g\in\RR[x]$ is called \emph{s.o.s.} polynomial if $g=\sum_{k=1}^N h_k^2$, for some $h_k\in\RR[x]$.
The subcone of s.o.s. polynomials in $\PHF{2d}{n}$ is denoted by $\SosHF{2d}{n}$;
it is a proper subcone  whenever $d\geq 2$, $n\geq 3$, unless $(d,n)=(2,3)$.
We use notation $f\geq 0$ to indicate global nonnegativity. Respectively, for inhomogeneous
$f$ we write $f>0$ to indicate global positivity, and for homogeneous $f$ we write $f>0$ to
indicate positivity away from the origin.

It is natural to consider a closely
related problem of finding a \emph{multiplier} $q\in \SosHF{2m}{n}$ for $f\in\PHF{2d}{n}$, so that
$qf\in\SosHF{2(d+m)}{n}$. Indeed, we can use the fact that the product of a
sum of squares is a sum of squares, and represent $f$ as a s.o.s. of rational
functions with common denominator $q$. Finding a multiplier may be cast as a \emph{s.o.s. polynomial
optimisation} feasibility problem, which is a case of \emph{semidefinite
optimisation} problem (SDP), cf. \cite[Sect.~3.1.7]{MR3075433},
along the lines of \cite{dKP04}. For the sake of completeness, we include a sketch of the corresponding
setup in Appendix~\ref{app:A}. In practice, these SDPs are (usually)
efficiently solvable using interior point methods, see e.g.  
\cite{MR3075433,MR3469431}---this was an important reason for a lot of attention to s.o.s. representations
in the past 25 years.

We show the following.
\begin{thm}\label{thm84}
Let $f\in\PHF{4}{4}$. Then $f=\frac{p}{qq'}$, 
where $p\in\SosHF{8}{4}$ and $q,q'\in\SosHF{2}{4}$.
\end{thm}

Note that while $q$, $q'$ might have some nontrivial real zeros,
the identity $f=\frac{p}{qq'}$ will hold almost everywhere in the usual measure-theoretic sense.

In a nutshell, the proof is as follows.
We start by explaining a folklore result that one can assume that $f$ has a real zero.
Then a linear change of variables allows reducing to the case where $f$ is a quadratic, w.r.t. $u:=x_1$
polynomial $f=au^2+2bu+c$. This means that for any fixed values of the
other variables (hidden in forms $a$, $b$, $c$), $f$ is a nonnegative quadratic polynomial, and so it suffices to
check that its minimum is nonnegative. We ``complete the square''
i.e. write $af=(au+b)^2+ac-b^2$, and observe that the discriminant $\Delta:=ac-b^2\in\PHF{6}{3}$.
As was shown by Hilbert  in \cite{Hilb-ternsos}, for any $h\in\PHF{6}{3}$
there exists a $q\in\SosHF{2}{3}$ so that $qh\in\SosHF{8}{3}$, thus such a $q$ exists for $\Delta$. Hence,
$aqf\in\SosHF{8}{4}$, as claimed. The details are in Section~\ref{sect:red0}.

In the remainder of the paper we investigate the sharpness of the degree
bound
in \Cref{thm84}.  Specifically,
as $\Delta$ lies in a proper subvariety of $\PHF{6}{3}$, for which no examples of non-s.o.s. nonnegative
forms were known prior to this work,
\emph{a priori} $\Delta$ could have always been a s.o.s., rendering the multiplier $q$ obsolete.
In Section~\ref{sect:tern} we construct, with help of Blekherman's \cite{MR2904568}, explicitly, non-s.o.s. $\Delta$ , cf. Theorem~\ref{thm:exDelta}. 
The following is an immediate consequence of the latter.
\begin{thm}\label{thm:nonsossextic}
There exist  $ac-b^2\in\PHF{6}{3}\setminus \SosHF{6}{3}$,
with  $a,b,c\in\RR[x,y,z]$, of degrees $2$, $3$, $4$, respectively. \qed
\end{thm}
Note that this result has been claimed by the author in versions 2 and 3 of \cite{DP2022}---however without a valid proof. Now, to facilitate reading, Appendix \ref{app:code} below includes code to construct and verify an example.

While examples from \Cref{sect:tern} provide a supply of potential
examples to show that a quadratic multiplier is not always sufficient,
for all of them a quadratic multiplier appears to be sufficient,
see the concluding \Cref{subsect:remarks}. In particular,
we describe there the relationship between s.o.s.
decompositions (with multipliers) of $\Delta$ and these of $f$. Namely,
\begin{thm}\label{thm:sos34}
Let $f:=au^2+2bu+c\in\PHF{4}{4}$ and $\Delta:=ac-b^2$. 
Then the existence of $q\in\SosHF{2}{3}$ so that $q^j f\in\SosHF{4+2j}{4},$
with $j\in\{0,1\}$, implies
$aq^j\Delta=\sum_{k}(aq_k-bq_k)^2\in\SosHF{8+2j}{3}$. 

Conversely, for $j=0$ the existence of the latter s.o.s. decomposition of $a\Delta$ implies that 
any $f=au^2+2b'u+c'$, with the discriminant equal to $\Delta$, is an s.o.s., i.e. $f\in\SosHF{4}{4}$.
\end{thm}
Note that the summands in the s.o.s. decomposition of
$aq^j\Delta$ in the above theorem belong to the ideal $(a,b)$.
We currently do not know how to prove the converse in the case $j=1$,
even in probably the most interesting case $q=a$---the latter would
be on a path to resolve a natural conjecture that $af\in\SosH{4}{6}$
always holds.

The condition $a\Delta\in (a,b)$ is necessary---examples of $\Delta=\Delta_{\tau}$ in 
\Cref{thm:exDelta} satisfy $a\Delta\in\SosHF{8}{3}$ for a range of values of $\tau$, but the corresponding $f$ is not in $\SosHF{4}{4}$.

Another such example, 
$\SosH{4}{4}\not\ni Q(u,y,z,x)
=(y^2+z^2)u^2-4xyzu+x^4+y^2z^2$,
so that $\Delta=(y^2+z^2)(x^4+y^2z^2)-(2xyz)^2=(y^2z-x^2z)^2+(x^2y-yz^2)^2$,
is provided in Choi and Lam \cite{ChLa77}.\footnote{In [loc.cit.] variables are named differently, $Q(x,y,z,w)=w^4+x^2y^2+y^2z^2+z^2x^2-4xyzw$.}

\section{Reduction to a real zero case}\label{sect:red0}
We start by showing that one can assume that $f$ has a real zero.
First, we recall a well-known lemma (e.g. similar observations
are made in \cite{MR2975385}).

\begin{lem}\label{lem1}
Let $f\in\PHF{2d}{n}$, 
and $S:=\{v\in \RR^{n}\mid \sum_{k=1}^n v_k^2=1\}$ be the unit sphere in $\RR^{n}$. Denote $u:=x_1$.
\begin{itemize}
\item[(i)] Let $f$ have a zero on $S$.  An orthogonal change
of coordinates bringing the zero to $(1{:}0{:}\dots {:} 0)$ transforms $f$ into the form
\begin{equation}
f(u,x)=\sum_{k=2}^{2d} f_k(x) u^{2d-k}, \quad f_k\in\RR[x_2,\dots,x_n].\label{red0}
\end{equation}
\item[(ii)] Let $f>0$, with a  minimum on $S$
reached at $x^*\in S$. 
An orthogonal change of coordinates bringing $x^*$ to $(1{:}0{:}\dots {:}0)$ transforms $f$ into the form
\begin{equation}
    f(u,x)=u^{2d}+\sum_{k=2}^{2d} f_k(x) u^{2d-k}, \quad f_k\in\RR[x_2,\dots,_n],\text{ so that $f(u,x)-u^{2d}\geq 0$.}\label{rednon0}
\end{equation}
\item[(iii)] $f_2\geq 0$ and $f_{2d}\geq 0$, and moreover, $f_2\in \SosHF{2}{n}$.
\end{itemize}
\end{lem}

\begin{proof}
(i). After the transformation, $f$ cannot have a term $u^{2d}$, as otherwise
it cannot vanish on $(1{:}0{:}\dots {:}0)$. It cannot have a term $f_1(x)u^{2d-1}$, as 
it is nonnegative. Thus we have the required form.

(ii).  Without loss of generality, after the transformation we have
\[
f(u,x)=u^{2d}+g(u,x), \quad g(u,x):=\sum_{k=1}^{2d} f_k(x) u^{2d-k}.
\]
Note that on $S$ one has $f(u,x)\geq 1$, as the minimum of $f$ on $S$ equals 1.
Thus on $S$ one has $g(u,x)=f(u,x)-u^{2d}\geq f(u,x)-1\geq 0$. Hence $g(u,x)\geq 0$,
as it is homogeneous. Moreover, it has a zero, $(1{:}0\dots {:}0)$, on 
$S$, and this (i) applies, ensuring $f_1=0$, as required.

(iii). $f_2\geq 0$ and $f_{2d}\geq 0$ follows from $f\geq 0$ in the case (i), and
from $f(u,x)-u^{2d}\geq 0$ in the case (ii). As $f_2$ is quadratic, it is an s.o.s.
\end{proof}

Therefore we may assume that $f$ has a real zero.
From now on, consider the case $d=2$. By Lemma \ref{lem1}, we have
\[ g(u,x):=au^2+2bu+c\geq 0,\quad a:=f_2(x))\geq 0,\ c:=f_4(x)\geq 0,\ b:=f_3(x)/2.\]
Note that the nonnegativity of $g$ implies that $ac -b^2\geq 0$.
Indeed, if $(x_2{:}\dots {:}x_n)$ were a point where $ac -b^2<0$, then the univariate
polynomial $g(u)=g(u,x_2,\dots,x_n)$ would have two real roots,
and would be negative on the interval between them, which is not possible.
Write, omitting $x$ for brevity:
\begin{equation}
ag(u)=a^2u^2+2abu+ac= (au+b)^2-b^2+ac. \label{abcsq}
\end{equation}
The latter is the sum of a square and a
polynomial $\Delta(x):=ac-b^2\in\PHF{6}{n}$.

In particular, for $n=3$, by \cite{Hilb-ternsos} (or by a recent \cite[Sect.
5-6]{MR3272733}, cf. also \cite{BSV19}), $\Delta(x)=\frac{u(x))}{q(x)}$,
with $u$ and $q$ s.o.s.  polynomials of degrees $8$ and $2$, respectively. To
summarise, we state

\begin{lem} \label{lem2} 
Let $g(u,x):=a(x)u^2+2b(x)u+c(x)$ be a nonnegative
degree $4$ homogeneous polynomial in $u,x=(x_2,x_3,x_4)$.
Then there exists $q(x)\in\SosHF{2}{3}$ such that
\begin{equation}
q(x)a(x)g(u,x)=\sum_{k=1}^N r_k(u,x)^2,\quad r_k\in \RR[u,x],\ \deg{r_k}=4, \label{soswith0}     
\end{equation}
i.e.  it is  an s.o.s. polynomial. \hfill\qed
\end{lem}

As an orthogonal change of coordinates (e.g. the inverse $G$ of
the one in Lemma \ref{lem1}(i)) respects s.o.s. 
decompositions, and as $g(1,0,0,0)=0$, we obtain, 
applying $G$ to the both sizes of \eqref{soswith0}, the following.

\begin{lem}\label{lem3}
Let $f\in\PHF{4}{4}$ have a zero on $S$. 
Then there exists 
$q_t(u,x)\in\SosHF{2}{4}$, $t=1,2$, such that
$q_1 q_2 f=s\in\SosHF{8}{4}$, and
$f=s/(q_1 q_2)$.  \hfill\qed
\end{lem}

To complete the proof of Theorem \ref{thm84}, 
it remains to observe that in the case of 
strictly positive $f(u,x)=u^4+f_2 u^2+f_3 u+f_4$ in
the form
\eqref{rednon0}, by Lemma \ref{lem2} and Lemma \ref{lem1}(ii) we have 
$f(u)=u^4+\frac{\sum_{k=1}^N r_k^2}{q f_2}$, i.e.

\[ q(u,x)f_2(u,x)f(u,x)=q(u,x)f_2(u,x)u^4+\sum_{k=1}^N r_k(u,x)^2.\]
The 1st term on the RHS of the latter is an s.o.s., and the argument
used to establish Lemma \ref{lem3}  applies here, too. \hfill\qed

\section{On sharpness of the s.o.s. degree bound for 4-ary 4-ics}\label{sect:tern}
Here we consider $f(x,y,z,u):=au^2+2bu+c$, where $a,b,c\in\RR[x,y,z]$, of degrees 2, 3,
4 respectively.
As is readily seen from \eqref{abcsq}, if $ac-b^2\in\SosHF{6}{3}$ then
it suffices to multiply $f$ by $a$ to obtain $af\in\SosHF{6}{4}$. Therefore
it is necessary, for non-existence of a quadratic multiplier to $f$
making $af$ s.o.s., that $ac-b^2\in\PHF{6}{3}\setminus\SosHF{6}{3}$.
Whether such sextics exist is a nontrivial question, not answered in the literature.
We answer it here in the affirmative. First, we settle the case of
the reducible discriminant.

\subsection{Reducible discriminant}\label{rem:easycases}
If $a=\ell^2$ then $f\in\SosHF{4}{4}$.
Indeed, as $\ell^2 c\geq b^2$, one has that $\ell$ divides $b$.
Thus one can write $b=b_1\ell$ and write down
\begin{equation*}
f=\ell^2 u^2+2b_1\ell u+c=(\ell u +b_1)^2+c-b_1^2,
\end{equation*}
As $\PHF{4}{3}=\SosHF{4}{3}$, $c-b_1^2$ is a s.o.s. Hence $f$ is a s.o.s. too.
(cf. \cite[Sect. 2.2]{MR2975385}).

Similarly, if $b=\ell a$, with $\ell$ linear, then
$f\in\SosHF{4}{4}$. Indeed,
\begin{equation*}
f=a u^2+2a\ell u+c=a(u +\ell)^2+c-a\ell^2,
\end{equation*}
and 
$c-a\ell^2$ is s.o.s., implying same for $f$.

Finally, let $c$ and $b^2$ have a nontrivial common factor, $q$, it must
have degree $2$ and satisfy $q\geq 0$, implying that $ac-b^2$ is an s.o.s.
implying that $af$ is a s.o.s., too.

\subsection{Preliminaries}
We assume, without loss of generality, that none of the cases in Section~\ref{rem:easycases} occur.
Not every polynomial in $\PHF{6}{3}$ can be written in the form $ac-b^2$, as can be observed
by the following dimension-counting argument, or extracted from a sheaf-theoretic argument from \cite{MR1786479}. Write
$ac-b^2=\left|\begin{matrix}a&b\cr b&c\end{matrix}\right|$ and observe that for any 3-ary linear form $\ell$
and a nonzero $r\in\RR$ one has
\begin{equation}\label{eq:cong}
\left|\begin{matrix}a&b\\ b&c\end{matrix}\right|=
\left|U^{\top}
\begin{pmatrix}a&b\\ b&c\end{pmatrix}
U\right|=
\left|\begin{matrix}r^2 a& r\ell a+b\\ r\ell a+b&\ell^2 a+2\ell b/r+c/r^2\end{matrix}\right|,\quad U:=\begin{pmatrix}r&\ell\\ 0&1/r\end{pmatrix}.
\end{equation}
Thus the variety $\mathcal{D}_-$ of nonnegative 3-ary sextics of the form $ac-b^2$ has dimension at most
$6+10+15-3-1=27<28$, where 28 is the dimension of $\PHF{6}{3}$.
Geometrically, a 3-ary sextic of the form $ac\pm b^2\in\mathcal{D}=\mathcal{D}_+\cup \mathcal{D}_-$ corresponds to a plane projective curve with double points of tangency to the conic $\{a=0\}$, and $\mathcal{D}, \mathcal{D}_+, \mathcal{D}_-$ have the same dimension, 27.

Note that all the representations of $\Delta:=ac-b^2$ with $\{a=0\}$ specifying the conic with 
the 6 double contact points  of the curve $\{\Delta=0\}$ given by the equations $a=b=0$  are essentially listed in \eqref{eq:cong}. In other words,
\begin{lem}\label{lem:cong}
Let the irreducible curve $\{\Delta=0\}$ have double contact points $\mathcal{C}:=\{C_1,\dots,C_6\}$ given by $a=b=0$. Then any representations of $\Delta$ of the form $a'c'-b'^2$,
with $a'$ vanishing on $\mathcal{C}$, is obtained from the action of an $U$.
\end{lem}
\begin{proof} There is unique conic through the contact points. Hence $a'=r^2a$ for some nonzero $r$. As $b'$ vanishes on $\mathcal{C}$, $b'=r\ell a+b$ for a linear form $\ell$. Now, $ac-b^2=r^2ac'-(r\ell a+b)^2$, implying
$a(c-c')=-r^2\ell^2 a^2-2r\ell a b$.
\end{proof}

At this point we are able to explain how $f$'s are related to their
discriminants more precisely. Namely, observe that 
\newcommand{\uv}{\begin{pmatrix} u\\ 1\end{pmatrix}}
$f(u):=au^2+2bu+c=\uv^\top \begin{pmatrix}a&b\\ b&c\end{pmatrix} \uv$
and $f_U(u):=\uv^\top U^\top \begin{pmatrix}a&b\\ b&c\end{pmatrix}U \uv$
have the same discriminant $\Delta$. On the other hand, let $f(u)$ and $g(u)=a'u^2+2b' u+c'$ have the same discriminant $\Delta$. 
If we know in addition that $a'=r^2a$ for some $r\in\RR$, then it follows
from Lemma~\ref{lem:cong} that $g(u)=f_U(u)$ for some $U$.
\begin{lem}\label{lem:cong44}
  Let $f(u)=au^2+2bu+c$ and $g(u)=r^2au^2+2b'u+c'$ satisfy $\Delta:=ac-b^2=r^2ac'-b'^2$. Then $g(u)=f_U(g)$. Also, $f_U(u)$ and $f(u)$ have
  the same $\Delta$. \qed
\end{lem}
It might, potentially, happen that $a'=b'=0$ give a different
from $\mathcal{C}$ set of double contact points of $\Delta$, and then we
cannot obtain $g$ as $f_U$. It is claimed in \cite[Chapter~I.3]{zbMATH06551558} that in this case $\Delta$ has infinitely many
sets of double contact points. Apparently, [loc.cit.]
is incorrect, see \cite{abx507635}.
As well, we do not know if a strictly positive $\Delta\in\partial\SosH{3}{6}$ with more than one
set of double contact points exists---it might be relevant here that the s.o.s. representation of strictly positive $\Delta\in\partial\SosH{3}{6}$ is essentially unique, see e.g. Capco and Scheiderer \cite{MR4275922}.

\subsection{Construction}
Returning to our plan to construct a non-s.o.s. $\Delta$, we proceed to find elements $p=ac-q_3^2$ in the intersection $\mathcal{D}_-$ with $\partial\SosH{3}{6}\setminus\partial\PH{3}{6}$, i.e. with the variety of extreme strictly positive s.o.s.,
as described in  \cite{MR2904568}, and further studied in \cite{algbound}.
With such $p$ at hand, $p-\epsilon q_3^2\in \PH{3}{6}\setminus \SosH{3}{6}$
for 
any sufficiently small $\epsilon>0$. 

On the other hand, any $p\in\partial\SosH{3}{6}\setminus\partial\PH{3}{6}$ is the sum $p=q_0^2+q_1^2+q_2^2$ of exactly 3 squares,  
cf. \cite{MR2904568}.
Note that $p_\zeta=\zeta_0q_0^2+\zeta_1q_1^2+\zeta_2q_2^2\in \partial\SosH{3}{6}\setminus\partial\PH{3}{6}$, as well, for any $\zeta>0$. 
Thus, we can hope to find $p_\zeta\in\mathcal{D}_-$, once we found $p$ in $\partial\SosH{3}{6}\setminus\partial\PH{3}{6}$, although we do not know whether this works for any such $p$.
However, it works for a specific example below. 
Namely, we show the following.
\begin{thm}\label{thm:exDelta}
Let $a=x^2+y^2+z^2$, 
$c=\frac{1}{2^4}\left(9x^4-5x^2y^2+36y^4-5x^2z^2-68y^2z^2+36z^4\right)$, $b=x(y^2-z^2)$.

Then
for any $4.6875=\frac{5^2 3}{2^4}<\tau\leq\tau^*=4.92231255059\dots$ one has $\Delta=\Delta_{\tau}:=ac-\tau b^2\in \PHF{6}{3}\setminus \SosHF{6}{3}$,
where $\tau^*$ is the unique real root of $m_t=36864t^3-213248^2+317872t-794375$,
while $0<\Delta_\frac{5^2 3}{2^4}\in \partial \SosH{3}{6}$.
\end{thm}

Consider the following, parametrised by $\alpha$, symmetric polynomial
\begin{equation}\label{eq:robinson}
  R_{\alpha}(x,y,z)=\alpha (x^6+y^6+z^6)-(x^2y^4+x^4y^2+x^2z^4+x^4z^2+y^2z^4+y^4z^2)+3x^2y^2z^2\in \PH{3}{6}, \quad \alpha\geq 1.
\end{equation}
Note that $R_1\in\partial\PH{3}{6}\setminus\SosH{3}{6}$ is a well-known example due to A.~Robinson,
cf. \cite{Rez-concr-h17}, with 10 real zeros $\{(1{:}\pm 1{:}\pm 1), (0{:}1{:}\pm 1), (1{:}0{:}\pm 1), (1{:}\pm 1{:}0)\}$.

\begin{lem}\label{lem:robinson}
  $2R_{\frac{9}{8}}=x^2(\frac{3}{2}x^2-y^2-z^2)^2+y^2(\frac{3}{2}y^2-x^2-z^2)^2+z^2(\frac{3}{2}z^2-x^2-y^2)^2=:q_0^2+q_1^2+q_2^2\in\partial\SosHF{3}{6}.$
\end{lem}
\begin{proof}
We demonstrate that $R_{\frac{9}{8}}\in\partial\SosH{3}{6}$ using a criterion from \cite{MR2904568}. Namely, consider the $9$ common roots
\begin{equation}\label{eq:q0q1roots}
Z:=\{z_1,\dots,z_9\}=\{(0:0:-1),(0:\sqrt{2}:\pm\sqrt{{3}}), (\sqrt{2}:0:\pm\sqrt{{3}}),
         (2:\pm 2:\pm {\sqrt{2}})\}
\end{equation}
of $q_0$ and $q_1$, given in \eqref{eq:q0q1roots}. There must exist a linear functional
$\ell:\PH{3}{6}\to \RR$ s.t. $\ell(\SosH{3}{6})\subseteq\RR_{\geq 0}$, $\ell(q_k^2)=0$ for $k=0,1,2$ given by $\ell(h)=\sum_{w\in Z}\mu_w h(w)$, with exactly one $\mu_w<0$, and the remaining eight $\mu_w>0$. Note that $\ell$ defines a quadric form $Q_{\ell}(q)=\sum_{w\in Z}\mu_w q(w)^2$, cf. \cite[Theorem~6.1]{MR2904568}.
The $\mu_w$'s can be computed from the values of $q_2$ on $Z$, and the 
\emph{Cayley-Bacharach coefficients} $u_w$, for $w\in Z$, that is, the coefficients
of the unique linear relation
\begin{multline}\label{eq:CB}
   \sum_{w\in Z}u_w q(w)=0, \quad 
   q\in H_{3,3}\text{ (with $H_{3,3}$ the $\mathbb{C}$-vectorspace of ternary cubic forms),}\\ (u_w)_{w\in Z}=(15\sqrt{2}, 2\sqrt{3}, -2\sqrt{3}, 2\sqrt{3}, -2\sqrt{3}, -3, 3, -3, 3),\quad
(\mu_w)_{w\in Z}=(-25, 4, 4, 4, 4, 6, 6, 6, 6)
\end{multline}
as $\mu_w=\frac{u_w}{q_2(w)}$, cf. \cite[Lemma~4.1]{MR2904568}, while
$(u_w)$ in \eqref{eq:CB} is obtained as the kernel of the $10\times 9$ matrix
of evaluations of a vectorspace basis elements of $H_{3,3}$ on $Z$. 
\end{proof}

\begin{rem}
The value $\alpha=\frac{9}{8}$ in Lemma~\ref{lem:robinson} was obtained by using a very convenient implementation \cite{SumsOfSquaresArticle,SumsOfSquaresSource} in Macaulay2 \cite{M2} of an algorithm to compute s.o.s. decompositions with
rational coefficients by Parrilo and Peyrl \cite{MR2474341}, by minimising the value of $\alpha$ for which $R_\alpha\in\SosH{3}{6}$. 

Results of computations in \eqref{eq:CB} were verified using \cite{sagemath}. 
\end{rem}

At this point, we have $p:=2R_{\frac{9}{8}}=p_{(1:1:1)}$, where
$p_\zeta=\zeta_0q_0^2+\zeta_1q_1^2+\zeta_2q_2^2$. As explained above, it remains to find
a value of $\zeta$ so that $p_\zeta\in\mathcal{D}_-$.
To this end, we restrict $p_\zeta$ to $\{a=0\}$ and seek $q_3$ s.t.  holds.
\begin{equation}\label{eq:onconic}
  p_\zeta+q_3^2=|_{a=0}0, \quad \zeta>0
\end{equation}
W.l.o.g. $q_3$ is a product of linear forms, allowing one to lift $q_3$ on $\{a=0\}$
to $H_{3,3}$.
As all variables are squared in $p$, restricting $p_\zeta$ to $\{a=0\}$ amounts to
substituting $z^2$ with $-x^2-y^2$.
A direct case analysis shows that up to scaling there are precisely 4 possibilities for
\eqref{eq:onconic} to hold. Namely, either all $\zeta_k=1$, or two of them equal to $1$, while the
3rd is equal to $\frac{1}{4}$. The case $\zeta=(1:1:1)$ corresponds to $q_3=Cxyz$, for some $C\in\RR$ and $p_{\zeta}\in\mathcal{D}_+$, not good for us. But $\zeta=(\frac{1}{4}:1:1)$ corresponds to
$q_3=\frac{75}{16}x(y^2-z^2)=\frac{75}{16}b$, and $p_{\zeta}\in\mathcal{D}_-$, as needed in 
the lower bound on $\tau$ in Theorem~\ref{thm:exDelta}. 
\Cref{fig:test} in \Cref{app:code} shows Macaulay2 code verifying that $p_\zeta\in\mathcal{D}_-$, with $a$, $b$, $c$, and $\tau=\frac{75}{16}$ as in \Cref{thm:exDelta}.

To get the upper bound on $\tau$, we consider the minima of the function
$p(x,y,z,t)=ac-tb^2$ on $\{a=1\}$. Namely, we look for the points where the gradients
of $p$ and of $a$ are collinear; we eliminate $x$, $y$, $z$, $\lambda$ from the ideal 
$J=(\frac{\partial p}{\partial x}-\lambda x, \frac{\partial p}{\partial y}-\lambda y, \frac{\partial p}{\partial z}-\lambda z, p, a-1)$,
after saturating $J$ with $x$ (as $x=0$ does not correspond to interesting values of $t$,
due to $b(0,y,z)=0$). The result, computed in Macaulay2 \cite{M2}\footnote{
Msolve \cite{berthomieu} used via its Macaulay2 interface \cite{MsolveSource}, was very helpful.}, the intersection of $I$ and $\RR[t]$, is generated by 
a degree 5 polynomial, the product of the minimal polynomial $m_t$ for $\tau^*$
and a strictly positive quadratic polynomial. This completes the proof of Theorem~\ref{thm:exDelta}. \qed

Using Macaulay2, we can find the real zeros of $\Delta_{\tau^*}$ as
the two components of the ideal $(\frac{\partial \Delta_{\tau^*}}{\partial x}, \frac{\partial \Delta_{\tau^*}}{\partial y}, \frac{\partial \Delta_{\tau^*}}{\partial z}, x-1)$, one of them
$(z,x-1,177003y^2+1152t^2-56920t+103528,m_t)$, and the other where
the roles of $y$ and $z$ are interchanged (note that $\Delta_{\tau}$ is invariant under swapping $y$ and $z$). Thus, the real projective zeros of $\Delta_{\tau^*}$ are as follows.
\begin{equation}\label{eq:deltarz}
  V_{\RR}(\Delta_{\tau^*})=\{(1:\pm\rho:0), (1:0:\pm\rho)\}, \quad 177003\rho^2+1152t^2-56920t+103528=0,\ m_t=0.
\end{equation}
A shorter representation for $\rho$, its minimal polynomial, is
\begin{equation} \label{eq:rhomin}
  18\rho^6-2\rho^2-9=0,\quad\text{giving } \tau^*=\frac{1}{144}(36\rho^4+486\rho^2+275).
\end{equation}

\section{Towards s.o.s. for quaternary quartics}\label{subsect:remarks}
The 4-ary 4-ic $f_*:=au^2+2\sqrt{\tau^*}bu+c$, with $a$, $b$, $c$, $\tau^*$ as in Theorem~\ref{thm:exDelta}
has 5 real projective roots\footnote{
According to \cite[Chapter~I.2]{zbMATH06551558}, any such quartic can be represented as a 5-ary quadratic form in 5 linearly independent quadrics through the 5 roots.}, namely $(0{:}0{:}0{:}1)$, and the 4 roots given by $\zeta\in V_{\RR}(\Delta_{\tau^*})$ in \eqref{eq:deltarz}, with the $u$-coordinate given by $-\frac{b(\zeta)}{a(\zeta)}$, as shown in \eqref{eq:extroots}, renormalised to
clean denominators.
However, $af_*$ appears to be an s.o.s., subject to the usual disclaimers related to the numerical computations done by Macaulay2 code in \Cref{app:code}, \Cref{fig:af_sos}.

We know for certain that $af(t)=a(au^2+2tbu+c)$ is an s.o.s., with rational coefficients,
for $t=\hat{t}:=\frac{2382234243}{1073741824}$, with $0<\tau^*-\hat{t}^2<3.10^{-9}$. This computation was done by maximising the value of $t$ subject to $af(t)$ being s.o.s., with subsequent rounding
of the answer to the exact rational values.
The question whether our degree bound is sharp remains open, for the time being,
as $ac-(\hat{t}b)^2$ is non-s.o.s. by Theorem~\ref{thm:exDelta}.

As $f_*$ has $5$ real roots, there are $5$ different choices for writing it out as a quadratic
polynomial with the corresponding discriminant; however, none of these discriminants turned
out to be s.o.s. Applying the change of coordinates given by the  matrix \eqref{eq:extroots} of the real projective
zeros of $\Delta_{\tau^*}$ 
\begin{equation} 
  \begin{pmatrix}\label{eq:extroots}
 1+ \rho^2 &  \rho(1+ \rho^2) &0 &- \tau^*  \rho^2\\
 1+ \rho^2 &- \rho (1+ \rho^2) &0 &- \tau^*  \rho^2\\
 1+ \rho^2 &0 & \rho (1+ \rho^2) & \tau^*  \rho^2\\
 1+ \rho^2 &0 &- \rho (1+ \rho^2) & \tau^*  \rho^2
 \end{pmatrix}
\end{equation}
to $f_*$ gives the polynomial invariant under the dihedral group of order 8
generated by involutions $(x,z)(y,u)$ and $(x,y)$. With the notation 
$f=\sum_{ijkl}f_{ijkl}x^iy^jz^ku^l$ and 

  \begin{align*}
\alpha&=(1192/9) \rho^4+(1168/9) \rho^2+89=f_{2200}=f_{0022}\\
\beta&=(14113/324) \rho^4+(1955/54) \rho^2+959/48=f_{0220}=f_{2020}=f_{2002}=f_{0202}\\
\gamma&=55 \rho^4+(89/2) \rho^2+27=f_{2110}=f_{1210}=f_{0121}=f_{2101}=f_{0112}=f_{1012}=f_{1201}=f_{1021}\\
\delta&=(4306/81) \rho^4+(6917/162) \rho^2+1855/72=f_{1120}=f_{2011}=f_{0211}=f_{1102}
\end{align*}
$f_*$ has the form
\begin{multline*}
f_*=\alpha(x^2y^2+z^2 u^2)+\beta(x^2 z^2+y^2 z^2+x^2 u^2+y^2 u^2)+
\gamma(x^2 y z+x y^2 z+x^2 y u+x y^2 u+x z^2 u+y z^2 u+x z u^2+y z u^2)\\
+\delta(x y z^2+x^2 z u+y^2 z u+x y u^2)
 -((17153/81) \rho^4+(17752/81) \rho^2+5219/36) x y z u. 
  \end{multline*}
These computations with $f_*$ may be verified with Macaulay2
code in \Cref{app:code}, \Cref{fig:transform}.

It might be helpful to obtain an exact s.o.s. decomposition for $af_*$ with as few terms as possible. This would need an extension of the rounding procedure from \cite{MR2474341}
to real number fields. One should also be aware of an
eventual need for doing approximate s.o.s. decompositions with
precision better than can be achieved by using standard SDP
solvers. 
One possible way out might be to use
an arbitrary precision SDP solver with \emph{controlled} precision, e.g. implementing an algorithm outlined in \cite{MR3549888}.

\begin{rem}
It is interesting to investigate a possibility to obtain a non-s.o.s. $ac-b^2\in\mathcal{D}_-$, with more than 4
real zeros (from  Theorem~\ref{thm:nonsossextic} is it not clear how to get more than 4 of them).
A theoretical bound is 10---this is how many zeroes ternary non-s.o.s. nonnegative sextics can have, see e.g.
\cite{algbound}. This bound is reached for $R_1\in\mathcal{D}_+$ \eqref{eq:robinson} (it can be checked that $R_1+Cx^2 y^2 z^2$, for some $C>0$, is divisible by $a$).

\end{rem}

\subsection{Proof of Theorem~\ref{thm:sos34}.} 
In view of \Cref{rem:easycases} one may assume that
$a$ does not divide $b$. As well, we may assume that the quadratic form $a$ is full rank---for the remaining cases the results follow by a continuity argument. We treat the case $j=0$ first.

Let $f=au^2+2bu+c=\frac{1}{a}((au+b)^2+\Delta)=\sum_k (p_ku+q_k)^2$, for $p_k$, $q_k$ forms of degree 1 and 2, respectively.
Plugging in $u=-\frac{b}{a}$, one obtains 
$\frac{\Delta}{a}=\sum_k (q_k-\frac{p_kb}{a})^2=\frac{1}{a^2}\sum_k (q_ka-p_kb)^2$, i.e. $a\Delta=\sum_k (q_ka-p_kb)^2$. That is, $f\in\SosH{4}{4}$ implies $a\Delta=\sum_k r_k^2\in\SosH{3}{8}$, with $r_k\in (a,b)$. Note that
$f=u^2\sum_k p_k^2+2u\sum_k p_k q_k+\sum_k q_k^2,$
i.e. $a=\sum_k p_k^2$, $b=\sum_k p_k q_k$, $c=\sum_k q_k^2$.

On the other hand, let $a\Delta=a^2c-ab^2=\sum_k r_k^2\in\SosH{3}{8}$, with $r_k=:q_k'a-p_k'b\in (a,b)$.

Then $a\Delta=\sum_k r_k^2=a^2\sum_{k}q_k'^2-2ab\sum_k p_k' q_k'+
b^2\sum_k p_k'^2$. As $b^2\sum_k p_k'^2$ must be divisible by $a$, we obtain
$\alpha a=\sum_k p_k'^2$, with $0\leq \alpha\in\RR$. As $\alpha=0$ implies that 
all the $p_k'=0$, and so $a$ divides $\Delta$, which is the case
where $\Delta\in\SosH{3}{6}$, and thus $af\in\SosH{4}{4}$ by \Cref{rem:easycases},
we may assume that $\alpha>0$. Dividing $a$ by $\alpha$ and at the same
time multiplying $c$ by $\alpha$ does not change $\Delta$, thus we may assume
that $\alpha=1$, i.e. $a=\sum_k p_k'^2$.

Also, note that $\frac{\Delta}{a}=\sum_k (q_k'-\frac{p_k'b}{a})^2$.
Substituting  $u=-\frac{b}{a}$ into the latter sum, we obtain
\begin{equation}\label{eq:fdelta}
f(u):=\sum_k(p_k'u+q_k')^2=au^2+2b'u+c', \quad\text{with }
a=\sum_k p_k'^2,\ b':=\sum_k p_k' q_k',\ c':=\sum_k q_k'^2.
\end{equation}
That is, there exists $f\in \SosH{4}{4}$ with given $a$ and $\Delta$.
By \Cref{lem:cong44}, any other such $g$ equals $f_U$ for
some $U$ in \eqref{eq:cong}. To complete the consideration of this case,
it remains to see that $f_U\in \SosH{4}{4}$ for any $U$.
Write 
\begin{multline*}
  \begin{pmatrix}a&b\\ b&c\end{pmatrix}=S^\top S,\quad f=\uv^\top S^\top S\uv,\qquad
  \text{where } S^\top:=\begin{pmatrix}p_1'&p_2'&\dots &p_m'\\
    q_1'&q_2'&\dots &q_m'\end{pmatrix},
    \quad U=\begin{pmatrix} r&\ell\\ 0&r^{-1}\end{pmatrix},
    \\
    \text{ and }\quad U^\top S^\top=\begin{pmatrix}rp_1'&rp_2'&\dots &rp_m'\\
    p_1'\ell+r^{-1}q_1'&p_1'\ell+r^{-1}q_1'&\dots &p_m'\ell+r^{-1}q_m'\end{pmatrix}.
\end{multline*}
Thus $f_U=\uv^\top U^\top S^\top SU\uv=\sum_{k=1}^m (rp'_k u +p_k'\ell+r^{-1}q_k')^2\in\SosH{4}{4}$,
and this case is done.

To complete the proof, let $qf=q(au^2+2bu+c)=\frac{q}{a}((au+b)^2+\Delta)=\sum_k (p_ku+q_k)^2$, for $p_k$, $q_k$ forms of degree 2 and 3, respectively.
Plugging in $u=-\frac{b}{a}$, one obtains 
$\frac{q\Delta}{a}=\sum_k (q_k-\frac{p_kb}{a})^2=\frac{1}{a^2}\sum_k (q_ka-p_kb)^2$, i.e. $aq\Delta=\sum_k (q_ka-p_kb)^2$. That is, $qf\in\SosH{4}{6}$ implies $aq\Delta=\sum_k r_k^2\in\SosH{3}{10}$, with $r_k\in (a,b)$. 
\qed

\subsection*{Acknowledgement} We thank Greg Blekherman, Alex Degtyarev, Claus Scheiderer, and Frank Sottile for helpful discussions.
The author was supported by the EU OpenDreamKit Horizon 2020 project and by ICERM, Brown University (funded by the National Science Foundation Grant No. DMS-1929284).

\bibliography{h17}
\bibliographystyle{abbrv}

\appendix
\section{An SDP for s.o.s. with a multiplier}\label{app:A}
Here we use multinomial notation for
monomials $X_1^{K_1}X_2^{K_2}\dots X_n^{K_n}:=X^K$
of degree $|K|:=\sum_j K_j$,
and the matrix scalar product $\langle A,B\rangle:=\Tr(A^\top B)=\sum\limits_{i,j}A_{ij}B_{ij}$.
Let $f(X)=\sum\limits_{|K|=2d}c_K X^k=\sum\limits_{t=1}^{t^*} (\sum\limits_{|M|=d} \ell_{t,M}X^M)^2\in\SosHF{2d}{n}$ be a sum of $t^*$ squares.
It is convenient to introduce a matrix $L$ by $L_{t,M}:=\ell_{t,M}$ and a vector $\mathcal{X}:=(X^M\mid |M|=d)$, allowing
one to write $f(X)=\mathcal{X}^\top L^\top L\mathcal{X}$. The matrix $G:=L^\top L$ is called a \emph{Gram matrix} for $f$.
By definition, $G$ is positive semidefinite (p.s.d, for short), and finding an s.o.s. decomposition for $f$ amounts to finding a $G$.
By comparing the coefficients of $X^K$ on the left- and the right-hand sides of $f(x)= \mathcal{X}^\top G\mathcal{X}$
we see that the latter is equivalent to the system of equations 
\begin{equation}\label{eq:gram}
\langle G,M_K\rangle =c_K,\qquad\text{ for each $K$ with $|K|=2d$},
\qquad
\end{equation}
where each $M_K$ is an explicit matrix which depends only on $K$; it is usually called \emph{moment matrix}.

Note that \eqref{eq:gram} together with the constraint saying that $G$ is p.s.d., written as $G\succeq 0$,
is a particular case of the feasibility problem for the semidefinite programming in primal form, cf. e.g. \cite[Chapter~2]{MR3075433}.

In a greater generality, where we only know that $f\geq 0$, we need to find
a  multiplier $a\in\SosHF{2m}{n}$, so that 
$af\in\SosHF{2(m+d)}{n}$. We put the Gram matrices for $af$ and $a$  into the block-diagonal p.s.d.
matrix variable $\diag(S,T)$. Now the right-hand sides of \eqref{eq:gram} become variables,
on which the condition that $a$ divides $af$ imposes linear constraints, and this is how the SDP
formulation for the problem at hand is usually stated in the literature. However, these extra variables can be
eliminated completely. Namely, the following holds.
\begin{prop}
Let $f\in\PHF{2d}{n}$. Then there exists $a\in\SosHF{2m}{n}$ satisfying
$af\in\SosHF{2(d+m)}{n}$, if and only if the linear system of equations 
\begin{equation}\label{eq:grammult}
\langle \diag(S,T), \diag(M_I,-\sum_{K+J=I}c_K M_J)\rangle = 0,\qquad\text{for $|I|=2(d+m)$},
\end{equation}
has a solution with $\Tr T=1$ and $T\succeq 0$, $S\succeq 0$.
\end{prop}
\begin{proof}
As $S$ is a Gram matrix for $af$, for every $I$ it must satisfy
\begin{equation*}\label{eq:grameq}
\langle S,M_I\rangle = \sum_{K+J=I} a_J c_K,\quad\text{where $a(X)=:\sum_J X^J$.}
\end{equation*}
On the other hand, $T$ is a Gram matrix for $a$, thus $a_J=\langle T,M_K\rangle$.
Therefore, by bi-linearity of $\langle ,\rangle$ we obtain
\begin{equation*}
\langle S,M_I\rangle = \sum_{K+J=I} c_K \langle T, M_J\rangle=\langle T,\sum_{K+J=I} c_K  M_J\rangle,
\end{equation*}
implying \eqref{eq:grammult}. The condition $\Tr T=1$ makes sure that $a$ is not identically 0.
\end{proof}

\clearpage

\section{Code to check computations}\label{app:code}
The code here is typeset using Macaulay2 package Merge\TeX \cite{MergeTeXSource}.

\begin{figure}[ht]
\begin{lstlisting}[language=Macaulay2]
R0 = frac QQ[r2,r3]/{r2^2-2, r3^2-3}
R = R0[x,y,z]
q0=x*((3/2)*x^2-y^2-z^2)
q1=y*((3/2)*y^2-x^2-z^2)
q2=z*((3/2)*z^2-x^2-y^2)
a=x^2+y^2+z^2
b=x*(y^2-z^2)
c=(1/2^4)*(9*x^4-5*x^2*y^2+36*y^4-5*x^2*z^2-68*y^2*z^2+36*z^4)
-- check that there is a contact conic {a=0}
assert(a*c-(75/16)*b^2==(1/4)*q0^2+q1^2+q2^2)
Z = apply({{0 , 0 , -1}, {0 , 1 , r3/r2}, {0 , 1 , -r3/r2}, {1 , 0 , r3/r2}, {1 , 0 , -r3/r2},
    {1 , 1 , r2/2}, {1 , 1 , -r2/2}, {1 , -1 , r2/2}, {1 , -1 , -r2/2}}, l->matrix(frac R0, {l}))
-- check that Z is the set of common roots of q0 and q1
assert(all(apply(Z, l->sub(q0,l)), m->m==0))
assert(all(apply(Z, l->sub(q1,l)), m->m==0))

-- check that q2 is non0 on each element of Z 
assert(all(apply(Z, l->sub(q2,l)), m->not m==0))

-- Cayley-Bacharach coefficients; check the relation on the monomial basis of qubics in R
u_w = {15*r2, 2*r3, -2*r3, 2*r3, -2*r3, -3, 3, -3, 3}
assert(all (apply ((entries basis(3,R))#0, nn-> sum(apply (u_w, Z, (m,l)->m*sub(nn,l)))), t->t==0))

-- quadratic form Q_l; check that it vanishes on q2^2
mu = {-25, 4, 4, 4, 4, 6, 6, 6, 6}
Q = q -> sum apply (mu, Z, (m,l)->m * (sub(q,l))^2)
assert(Q q2 == 0)
assert(Q b^2 > 0) -- but is positive on b^2
\end{lstlisting}
\caption{Computations for the proof of \Cref{thm:exDelta}.\label{fig:test}}
\end{figure}

\begin{figure}[ht]
\begin{lstlisting}[language=Macaulay2]
loadPackage "SumsOfSquares"
R = QQ[x,y,z,u][t]
a=x^2+y^2+z^2
b=x*(y^2-z^2)
c=(1/2^4)*(9*x^4-5*x^2*y^2+36*y^4-5*x^2*z^2-68*y^2*z^2+36*z^4)
F=a*u^2+2*t*b*u+c
s=solveSOS(F,-t)
assert(s#Parameters_0_0==55/32)
assert(value sosPoly s == sub(F,{t=>55/32})) -- exact check; for this value t we have s.o.s.
s=solveSOS(a*F,-t)
tt=s#Parameters_0_0
assert(tt==2382234243/1073741824) -- tt^2 ~ 4.922312547597944 (about 3.10^-9 smaller than \tau^*)
\end{lstlisting}
\caption{Computations to check that $af_*$ in \Cref{subsect:remarks} is s.o.s.
\label{fig:af_sos}}
\end{figure}

\begin{figure}[ht]
\begin{lstlisting}[language=Macaulay2]
RC=QQ[t,p, MonomialOrder => Eliminate 1]
Rc=RC/ideal(18*p^6-2*p^2-9, 144*t^2-36*p^4-486*p^2-275)
R=Rc[x,y,z]
a=x^2+y^2+z^2
b=x*(y^2-z^2)
c=(1/2^4)*(9*x^4-5*x^2*y^2+36*y^4-5*x^2*z^2-68*y^2*z^2+36*z^4)
f=a*c-t^2*b^2
m=matrix{{1,p,0},{1,-p,0},{1,0,p}}
su= (transpose m)*vector{x,y,z}
sub (f, transpose matrix{su})
Ru=Rc[x,y,z,u]
mRu=map(Ru,R)
a=mRu(a)
b=mRu(b)
c=mRu(c)
f=a*c-t^2*b^2
F=a*u^2+2*t*b*u+c
m=matrix{{1+p^2, p*(1+p^2),0,-t*p^2},
         {1+p^2,-p*(1+p^2),0,-t*p^2},
         {1+p^2,0,p*(1+p^2),t*p^2},
         {1+p^2,0,-p*(1+p^2),t*p^2}}
su= (transpose m)*vector{x,y,z,u}

rr = sub (F, transpose matrix{su})
assert( sub (rr, {x=>1, y=>1, z=>-1, u=>-1})==0) -- 5th root
assert( sub (rr, {x=>y, y=>x})==rr) -- D_8 - symmetric
assert( sub (rr, {u=>z, z=>u})==rr)
assert( sub (rr, {x=>z, y=>u, z=>x, u=>y})==rr)
\end{lstlisting}
\caption{Computations for the polynomial $f_*$ in \Cref{subsect:remarks}.
\label{fig:transform}}
\end{figure}

\end{document}

%% file: lst-Macaulay2.tex
\usepackage{listings}
\usepackage{xcolor}
\definecolor{maccolor}{rgb}{0.3,0.3,0.8}
\lstdefinelanguage{Macaulay2}
{
basicstyle={\ttfamily},
keywordstyle={\color{maccolor!80!black}},
commentstyle={\color{gray}},
stringstyle={\color{red!40!black}},
rulecolor=\color{maccolor},
basewidth={1.2ex}, 
sensitive=false,
morecomment=[l]{--},
morecomment=[s]{-*}{*-},
morestring=[b]",
escapechar={`},
escapebegin={\rmfamily},
morekeywords={about,abs,AbstractToricVarieties,accumulate,Acknowledgement,acos,acosh,acot,addCancelTask,addDependencyTask,addEndFunction,addHook,AdditionalPaths,addStartFunction,addStartTask,Adjacent,adjoint,AdjointIdeal,AffineVariety,AfterEval,AfterNoPrint,AfterPrint,agm,AInfinity,alarm,AlgebraicSplines,Algorithm,Alignment,all,AllCodimensions,allowableThreads,ambient,analyticSpread,Analyzer,AnalyzeSheafOnP1,ancestor,ancestors,ANCHOR,and,andP,AngleBarList,ann,annihilator,antipode,any,append,applicationDirectory,applicationDirectorySuffix,apply,applyKeys,applyPairs,applyTable,applyValues,apropos,argument,Array,arXiv,Ascending,ascii,asin,asinh,ass,assert,associatedGradedRing,associatedPrimes,AssociativeAlgebras,AssociativeExpression,atan,atan2,atEndOfFile,Authors,autoload,AuxiliaryFiles,backtrace,Bag,Bareiss,baseFilename,BaseFunction,baseName,baseRing,baseRings,BaseRow,BasicList,basis,BasisElementLimit,Bayer,BeforePrint,beginDocumentation,BeginningMacaulay2,Benchmark,benchmark,Bertini,BesselJ,BesselY,betti,BettiCharacters,BettiTally,between,BGG,BIBasis,Binary,BinaryOperation,Binomial,binomial,BinomialEdgeIdeals,Binomials,BKZ,BlockMatrix,BLOCKQUOTE,BODY,Body,BoijSoederberg,BOLD,Book3264Examples,Boolean,BooleanGB,borel,Boxes,BR,break,Browse,Bruns,cache,CacheExampleOutput,CacheFunction,CacheTable,cacheValue,CallLimit,cancelTask,capture,catch,Caveat,CC,CDATA,ceiling,Center,centerString,Certification,ChainComplex,chainComplex,ChainComplexExtras,ChainComplexMap,ChainComplexOperations,ChangeMatrix,char,CharacteristicClasses,characters,charAnalyzer,check,CheckDocumentation,chi,Chordal,class,Classic,clean,clearAll,clearEcho,clearOutput,close,closeIn,closeOut,ClosestFit,CODE,code,codim,CodimensionLimit,coefficient,CoefficientRing,coefficientRing,coefficients,Cofactor,CohenEngine,CohenTopLevel,CoherentSheaf,CohomCalg,cohomology,coimage,CoincidentRootLoci,coker,cokernel,collectGarbage,columnAdd,columnate,columnMult,columnPermute,columnRankProfile,columnSwap,combine,Command,commandInterpreter,commandLine,COMMENT,commonest,commonRing,comodule,CompactMatrix,compactMatrixForm,CompiledFunction,CompiledFunctionBody,CompiledFunctionClosure,Complement,complement,complete,CompleteIntersection,CompleteIntersectionResolutions,Complexes,ComplexField,components,compose,compositions,compress,concatenate,conductor,ConductorElement,cone,Configuration,ConformalBlocks,conjugate,connectionCount,Consequences,Constant,Constants,constParser,content,continue,contract,Contributors,ConvexInterface,conwayPolynomial,ConwayPolynomials,copy,copyDirectory,copyFile,copyright,Core,CorrespondenceScrolls,cos,cosh,cot,CotangentSchubert,cotangentSheaf,coth,cover,coverMap,cpuTime,createTask,Cremona,csc,csch,current,currentColumnNumber,currentDirectory,currentFileDirectory,currentFileName,currentLayout,currentLineNumber,currentPackage,currentString,currentTime,Cyclotomic,Database,Date,DD,dd,deadParser,debug,debugError,DebuggingMode,debuggingMode,debugLevel,DecomposableSparseSystems,Decompose,decompose,deepSplice,Default,default,defaultPrecision,Degree,degree,degreeLength,DegreeLift,DegreeLimit,DegreeMap,DegreeOrder,DegreeRank,Degrees,degrees,degreesMonoid,degreesRing,delete,demark,denominator,Dense,Density,Depth,depth,Descending,Descent,Describe,describe,Description,det,determinant,DeterminantalRepresentations,DGAlgebras,diagonalMatrix,diameter,Dictionary,dictionary,dictionaryPath,diff,DiffAlg,difference,dim,directSum,disassemble,discriminant,dismiss,Dispatch,distinguished,DIV,Divide,divideByVariable,DivideConquer,DividedPowers,Divisor,DL,Dmodules,do,doc,docExample,docTemplate,document,DocumentTag,Down,drop,DT,dual,eagonNorthcott,EagonResolution,echoOff,echoOn,EdgeIdeals,edit,EigenSolver,eigenvalues,eigenvectors,eint,EisenbudHunekeVasconcelos,elapsedTime,elapsedTiming,elements,Eliminate,eliminate,Elimination,EliminationMatrices,EllipticCurves,EllipticIntegrals,else,EM,Email,End,end,endl,endPackage,Engine,engineDebugLevel,EngineRing,EngineTests,entries,EnumerationCurves,environment,Equation,EquivariantGB,erase,erf,erfc,error,errorDepth,euler,EulerConstant,eulers,even,EXAMPLE,ExampleFiles,ExampleItem,examples,ExampleSystems,Exclude,exec,exit,exp,expectedReesIdeal,expm1,exponents,export,exportFrom,exportMutable,Expression,expression,Ext,extend,ExteriorIdeals,ExteriorModules,exteriorPower,Factor,factor,false,Fano,FastMinors,FastNonminimal,FGLM,File,fileDictionaries,fileExecutable,fileExists,fileExitHooks,fileLength,fileMode,FileName,FilePosition,fileReadable,fileTime,fileWritable,fillMatrix,findFiles,findHeft,FindOne,findProgram,findSynonyms,FiniteFittingIdeals,First,first,firstkey,FirstPackage,fittingIdeal,flagLookup,FlatMonoid,flatten,flattenRing,Flexible,flip,floor,flush,fold,FollowLinks,for,forceGB,fork,FormalGroupLaws,Format,format,formation,FourierMotzkin,FourTiTwo,fpLLL,frac,fraction,FractionField,frames,FrobeniusThresholds,from,fromDividedPowers,fromDual,Function,FunctionApplication,FunctionBody,functionBody,FunctionClosure,FunctionFieldDesingularization,fusePairs,futureParser,GaloisField,Gamma,gb,GBDegrees,gbRemove,gbSnapshot,gbTrace,gcd,gcdCoefficients,gcdLLL,GCstats,genera,GeneralOrderedMonoid,GenerateAssertions,generateAssertions,generator,generators,Generic,GenericInitialIdeal,genericMatrix,genericSkewMatrix,genericSymmetricMatrix,gens,genus,get,getc,getChangeMatrix,getenv,getGlobalSymbol,getNetFile,getNonUnit,getPrimeWithRootOfUnity,getSymbol,getWWW,GF,gfanInterface,Givens,GKMVarieties,GLex,Global,global,globalAssign,globalAssignFunction,GlobalAssignHook,globalAssignment,globalAssignmentHooks,GlobalDictionary,GlobalHookStore,globalReleaseFunction,GlobalReleaseHook,Gorenstein,GradedLieAlgebras,GradedModule,gradedModule,GradedModuleMap,gradedModuleMap,gramm,GraphicalModels,GraphicalModelsMLE,Graphics,graphIdeal,graphRing,Graphs,Grassmannian,GRevLex,GroebnerBasis,groebnerBasis,GroebnerBasisOptions,GroebnerStrata,GroebnerWalk,groupID,GroupLex,GroupRevLex,GTZ,Hadamard,handleInterrupts,HardDegreeLimit,hash,HashTable,hashTable,HEAD,HEADER1,HEADER2,HEADER3,HEADER4,HEADER5,HEADER6,HeaderType,Heading,Headline,Heft,heft,Height,height,help,Hermite,hermite,Hermitian,HH,hh,HigherCIOperators,HighestWeights,Hilbert,hilbertFunction,hilbertPolynomial,hilbertSeries,HodgeIntegrals,hold,Holder,Hom,homeDirectory,HomePage,Homogeneous,Homogeneous2,homogenize,homology,homomorphism,HomotopyLieAlgebra,hooks,horizontalJoin,HorizontalSpace,HR,HREF,HTML,html,httpHeaders,Hybrid,HyperplaneArrangements,Hypertext,hypertext,HypertextContainer,HypertextParagraph,icFracP,icFractions,icMap,icPIdeal,id,Ideal,ideal,idealizer,identity,if,IgnoreExampleErrors,ii,image,imaginaryPart,IMG,ImmutableType,importFrom,in,incomparable,Increment,independentSets,indeterminate,IndeterminateNumber,Index,index,indexComponents,IndexedVariable,IndexedVariableTable,indices,inducedMap,inducesWellDefinedMap,InexactField,InexactFieldFamily,InexactNumber,InfiniteNumber,infinity,info,InfoDirSection,infoHelp,Inhomogeneous,input,Inputs,insert,installAssignmentMethod,installedPackages,installHilbertFunction,installMethod,installMinprimes,installPackage,InstallPrefix,instance,instances,IntegralClosure,integralClosure,integrate,IntermediateMarkUpType,interpreterDepth,intersect,intersectInP,Intersection,intersection,interval,InvariantRing,inverse,InverseMethod,inversePermutation,Inverses,inverseSystem,InverseSystems,Invertible,InvolutiveBases,irreducibleCharacteristicSeries,irreducibleDecomposition,isAffineRing,isANumber,isBorel,isCanceled,isCommutative,isConstant,isDirectory,isDirectSum,isEmpty,isField,isFinite,isFinitePrimeField,isFreeModule,isGlobalSymbol,isHomogeneous,isIdeal,isInfinite,isInjective,isInputFile,isIsomorphism,isLinearType,isListener,isLLL,isMember,isModule,isMonomialIdeal,isNormal,isOpen,isOutputFile,isPolynomialRing,isPrimary,isPrime,isPrimitive,isPseudoprime,isQuotientModule,isQuotientOf,isQuotientRing,isReady,isReal,isReduction,isRegularFile,isRing,isSkewCommutative,isSorted,isSquareFree,isStandardGradedPolynomialRing,isSubmodule,isSubquotient,isSubset,isSupportedInZeroLocus,isSurjective,isTable,isUnit,isWellDefined,isWeylAlgebra,ITALIC,Iterate,Jacobian,jacobian,jacobianDual,Jets,Join,join,Jupyter,K3Carpets,K3Surfaces,Keep,KeepFiles,KeepZeroes,ker,kernel,kernelLLL,kernelOfLocalization,Key,keys,Keyword,Keywords,kill,koszul,Kronecker,KustinMiller,LABEL,last,lastMatch,LATER,LatticePolytopes,Layout,lcm,leadCoefficient,leadComponent,leadMonomial,leadTerm,Left,left,length,LengthLimit,letterParser,Lex,LexIdeals,LI,Licenses,LieTypes,lift,liftable,Limit,limitFiles,limitProcesses,Linear,LinearAlgebra,LinearTruncations,lineNumber,lines,LINK,linkFile,List,list,listForm,listLocalSymbols,listSymbols,listUserSymbols,LITERAL,LLL,LLLBases,lngamma,load,loadDepth,LoadDocumentation,loadedFiles,loadedPackages,loadPackage,Local,local,localDictionaries,LocalDictionary,localize,LocalRings,locate,log,log1p,LongPolynomial,lookup,lookupCount,LowerBound,LUdecomposition,M0nbar,M2CODE,Macaulay2Doc,makeDirectory,MakeDocumentation,makeDocumentTag,MakeHTML,MakeInfo,MakeLinks,makePackageIndex,MakePDF,makeS2,Manipulator,map,MapExpression,MapleInterface,markedGB,Markov,MarkUpType,match,mathML,Matrix,matrix,MatrixExpression,Matroids,max,maxAllowableThreads,maxExponent,MaximalRank,maxPosition,MaxReductionCount,MCMApproximations,member,memoize,memoizeClear,memoizeValues,MENU,merge,mergePairs,META,method,MethodFunction,MethodFunctionBinary,MethodFunctionSingle,MethodFunctionWithOptions,methodOptions,methods,midpoint,min,minExponent,mingens,mingle,minimalBetti,MinimalGenerators,MinimalMatrix,minimalPresentation,minimalPresentationMap,minimalPresentationMapInv,MinimalPrimes,minimalPrimes,minimalReduction,Minimize,minimizeFilename,MinimumVersion,minors,minPosition,minPres,minprimes,Minus,minus,Miura,MixedMultiplicity,mkdir,mod,Module,module,ModuleDeformations,modulo,MonodromySolver,Monoid,monoid,MonoidElement,Monomial,MonomialAlgebras,monomialCurveIdeal,MonomialIdeal,monomialIdeal,MonomialIntegerPrograms,MonomialOrbits,MonomialOrder,Monomials,monomials,MonomialSize,monomialSubideal,moveFile,multidegree,multidoc,multigraded,MultigradedBettiTally,MultiGradedRationalMap,multiplicity,MultiplicitySequence,MultiplierIdeals,MultiplierIdealsDim2,MultiprojectiveVarieties,mutable,MutableHashTable,mutableIdentity,MutableList,MutableMatrix,mutableMatrix,NAGtypes,Name,nanosleep,Nauty,NautyGraphs,NCAlgebra,NCLex,needs,needsPackage,Net,net,NetFile,netList,new,newClass,newCoordinateSystem,NewFromMethod,newline,NewMethod,newNetFile,NewOfFromMethod,NewOfMethod,newPackage,newRing,nextkey,nextPrime,nil,NNParser,NoetherianOperators,NoetherNormalization,NonminimalComplexes,nonspaceAnalyzer,NoPrint,norm,normalCone,Normaliz,NormalToricVarieties,not,Nothing,notify,notImplemented,NTL,null,nullaryMethods,nullhomotopy,nullParser,nullSpace,Number,number,NumberedVerticalList,numcols,numColumns,numerator,numeric,NumericalAlgebraicGeometry,NumericalCertification,NumericalImplicitization,NumericalLinearAlgebra,NumericalSchubertCalculus,numericInterval,NumericSolutions,numgens,numRows,numrows,odd,oeis,of,ofClass,OL,OldPolyhedra,OldToricVectorBundles,on,OneExpression,OnlineLookup,OO,oo,ooo,oooo,openDatabase,openDatabaseOut,openFiles,openIn,openInOut,openListener,OpenMath,openOut,openOutAppend,operatorAttributes,Option,OptionalComponentsPresent,optionalSignParser,Options,options,OptionTable,optP,or,Order,order,OrderedMonoid,orP,OutputDictionary,Outputs,override,pack,Package,package,PackageCitations,PackageDictionary,PackageExports,PackageImports,PackageTemplate,packageTemplate,pad,pager,PairLimit,pairs,PairsRemaining,PARA,Parametrization,parent,Parenthesize,Parser,Parsing,part,Partition,partition,partitions,parts,path,pdim,peek,PencilsOfQuadrics,Permanents,permanents,permutations,pfaffians,PHCpack,PhylogeneticTrees,pi,PieriMaps,pivots,PlaneCurveSingularities,plus,poincare,poincareN,Points,polarize,poly,Polyhedra,Polymake,PolynomialRing,Posets,Position,position,positions,PositivityToricBundles,POSIX,Postfix,Power,power,powermod,PRE,Precision,precision,Prefix,prefixDirectory,prefixPath,preimage,prepend,presentation,pretty,primaryComponent,PrimaryDecomposition,primaryDecomposition,PrimaryTag,PrimitiveElement,Print,print,printerr,printingAccuracy,printingLeadLimit,printingPrecision,printingSeparator,printingTimeLimit,printingTrailLimit,printString,printWidth,processID,Product,product,ProductOrder,profile,profileSummary,Program,programPaths,ProgramRun,Proj,Projective,ProjectiveHilbertPolynomial,projectiveHilbertPolynomial,ProjectiveVariety,promote,protect,Prune,prune,PruneComplex,pruningMap,Pseudocode,pseudocode,pseudoRemainder,Pullback,PushForward,pushForward,Python,QQ,QQParser,QRDecomposition,QthPower,Quasidegrees,QuaternaryQuartics,QuillenSuslin,quit,Quotient,quotient,quotientRemainder,QuotientRing,Radical,radical,RadicalCodim1,radicalContainment,RaiseError,random,RandomCanonicalCurves,RandomComplexes,RandomCurves,RandomCurvesOverVerySmallFiniteFields,RandomGenus14Curves,RandomIdeals,randomKRationalPoint,RandomMonomialIdeals,randomMutableMatrix,RandomObjects,RandomPlaneCurves,RandomPoints,RandomSpaceCurves,Range,rank,RationalMaps,RationalPoints,RationalPoints2,ReactionNetworks,read,readDirectory,readlink,readPackage,RealField,RealFP,realPart,realpath,RealQP,RealQP1,RealRoots,RealRR,RealXD,recursionDepth,recursionLimit,Reduce,reducedRowEchelonForm,reduceHilbert,reductionNumber,ReesAlgebra,reesAlgebra,reesAlgebraIdeal,reesIdeal,References,ReflexivePolytopesDB,regex,regexQuote,registerFinalizer,regSeqInIdeal,Regularity,regularity,relations,RelativeCanonicalResolution,relativizeFilename,Reload,remainder,RemakeAllDocumentation,remove,removeDirectory,removeFile,removeLowestDimension,reorganize,replace,RerunExamples,res,reshape,ResidualIntersections,ResLengthThree,Resolution,resolution,ResolutionsOfStanleyReisnerRings,restart,Result,resultant,Resultants,return,returnCode,Reverse,reverse,RevLex,Right,right,Ring,ring,RingElement,RingFamily,ringFromFractions,RingMap,rootPath,roots,rootURI,rotate,round,rowAdd,RowExpression,rowMult,rowPermute,rowRankProfile,rowSwap,RR,RRi,rsort,run,RunDirectory,RunExamples,RunExternalM2,runHooks,runLengthEncode,runProgram,same,saturate,Saturation,scan,scanKeys,scanLines,scanPairs,scanValues,schedule,schreyerOrder,Schubert,Schubert2,SchurComplexes,SchurFunctors,SchurRings,SCRIPT,scriptCommandLine,ScriptedFunctor,SCSCP,searchPath,sec,sech,SectionRing,SeeAlso,seeParsing,SegreClasses,select,selectInSubring,selectVariables,SelfInitializingType,SemidefiniteProgramming,Seminormalization,separate,SeparateExec,separateRegexp,Sequence,sequence,Serialization,serialNumber,Set,set,setEcho,setGroupID,setIOExclusive,setIOSynchronized,setIOUnSynchronized,setRandomSeed,setup,setupEmacs,sheaf,SheafExpression,sheafExt,sheafHom,SheafOfRings,shield,ShimoyamaYokoyama,short,show,showClassStructure,showHtml,showStructure,showTex,showUserStructure,SimpleDoc,simpleDocFrob,SimplicialComplexes,SimplicialDecomposability,SimplicialPosets,SimplifyFractions,sin,singularLocus,sinh,size,size2,SizeLimit,SkewCommutative,SlackIdeals,sleep,SLnEquivariantMatrices,SLPexpressions,SMALL,smithNormalForm,solve,someTerms,Sort,sort,sortColumns,SortStrategy,source,SourceCode,SourceRing,SPACE,SpaceCurves,SPAN,span,SparseMonomialVectorExpression,SparseResultants,SparseVectorExpression,Spec,SpechtModule,SpecialFanoFourfolds,specialFiber,specialFiberIdeal,SpectralSequences,splice,splitWWW,sqrt,SRdeformations,stack,stacksProject,Standard,standardForm,standardPairs,StartWithOneMinor,stashValue,StatePolytope,StatGraphs,status,stderr,stdio,step,StopBeforeComputation,stopIfError,StopWithMinimalGenerators,Strategy,String,STRONG,StronglyStableIdeals,STYLE,Style,style,SUB,sub,SubalgebraBases,sublists,submatrix,submatrixByDegrees,Subnodes,subquotient,SubringLimit,Subscript,subscript,SUBSECTION,subsets,substitute,substring,subtable,Sugarless,Sum,sum,SumOfTwists,SumsOfSquares,SUP,super,SuperLinearAlgebra,Superscript,superscript,support,SVD,SVDComplexes,switch,SwitchingFields,sylvesterMatrix,Symbol,symbol,SymbolBody,symbolBody,SymbolicPowers,symlinkDirectory,symlinkFile,symmetricAlgebra,symmetricAlgebraIdeal,symmetricKernel,SymmetricPolynomials,symmetricPower,synonym,SYNOPSIS,syz,Syzygies,SyzygyLimit,SyzygyMatrix,SyzygyRows,syzygyScheme,TABLE,Table,table,take,Tally,tally,tan,TangentCone,tangentCone,tangentSheaf,tanh,target,Task,taskResult,TateOnProducts,TD,temporaryFileName,tensor,tensorAssociativity,TensorComplexes,terminalParser,terms,TEST,Test,testExample,testHunekeQuestion,TestIdeals,TestInput,tests,TEX,tex,TeXmacs,texMath,Text,TH,then,Thing,ThinSincereQuivers,ThreadedGB,threadVariable,Threshold,throw,Time,time,times,timing,TITLE,TO,to,TO2,toAbsolutePath,toCC,toDividedPowers,toDual,toExternalString,toField,TOH,toList,toLower,top,top,topCoefficients,Topcom,topComponents,topLevelMode,Tor,TorAlgebra,Toric,ToricInvariants,ToricTopology,ToricVectorBundles,toRR,toRRi,toSequence,toString,TotalPairs,toUpper,TR,trace,transpose,TriangularSets,Tries,Trim,trim,Triplets,Tropical,true,Truncate,truncate,truncateOutput,Truncations,try,TSpreadIdeals,TT,tutorial,Type,TypicalValue,typicalValues,UL,ultimate,unbag,uncurry,Undo,undocumented,uniform,uninstallAllPackages,uninstallPackage,Unique,unique,Units,Unmixed,unsequence,unstack,Up,UpdateOnly,UpperTriangular,URL,urlEncode,Usage,use,UseCachedExampleOutput,UseHilbertFunction,UserMode,userSymbols,UseSyzygies,utf8,utf8check,validate,value,values,Variable,VariableBaseName,Variables,Variety,variety,vars,Vasconcelos,Vector,vector,VectorExpression,VectorFields,VectorGraphics,Verbose,Verbosity,Verify,VersalDeformations,versalEmbedding,Version,version,VerticalList,VerticalSpace,viewHelp,VirtualResolutions,VirtualTally,VisibleList,Visualize,wait,WebApp,wedgeProduct,weightRange,Weights,WeylAlgebra,WeylGroups,when,whichGm,while,width,wikipedia,Wrap,wrap,WrapperType,XML,xor,youngest,zero,ZeroExpression,zeta,ZZ,ZZParser}
}
\lstalias{Macaulay2output}{Macaulay2}